\documentclass{article}

\usepackage{amsmath}
\usepackage{amsfonts}
\usepackage{amsthm}
\usepackage{stmaryrd}
\usepackage{rotating}
\usepackage{pifont}
\usepackage{enumerate}

\theoremstyle{plain}
\newtheorem{theorem}{Theorem}[section]

\newtheorem{proposition}[theorem]{Proposition}
\newtheorem{corollary}[theorem]{Corollary}
\newtheorem*{conjecture}{Conjecture}

\theoremstyle{definition}
\newtheorem*{definition}{Definition}

\newcommand{\cd}{{\bf cd}}
\renewcommand{\c}{\ensuremath{\mathbf{c}}}
\renewcommand{\d}{\ensuremath{\mathbf{d}}}
\newcommand{\cyclic}[2]{\ensuremath{\mathcal{C}_{#1}(#2)}}

\newcommand{\polar}[1]{\ensuremath{{#1}^{\Delta}}}
\newcommand{\convexity}{{\rm (C)}}
\newcommand{\logconvexity}{{\rm (L)}}
\newcommand{\unimodality}{{\rm (U)}}
\newcommand{\barany}{{\rm (B)}}
\newcommand{\holds}{\ding{52}}
\newcommand{\holdsnot}{\ding{56}}

\begin{document}

\title{Unimodality and convexity \\ of $f$-vectors of polytopes}
\author{Axel Werner
  \thanks{TU Berlin, Institute of Mathematics, Berlin (Germany),
    {\tt awerner@math.tu-berlin.de}}}
\date{December 6, 2005}
\maketitle

\begin{abstract}
We consider unimodality and related properties of $f$-vectors of polytopes
in various dimensions. By a result of Kalai (1988), $f$-vectors
of $5$-polytopes are unimodal. In higher dimensions much less can be said;
we give an overview on current results and present a potentially interesting
construction as well as a conjecture arising from this.
\end{abstract}


\section{Introduction} \label{sec:introduction}

Let $f = (f_0,\ldots,f_{d-1})$ be the $f$-vector of a $d$-polytope.
It is natural to ask whether the $f$-vector necessarily has one (or more) of the following properties:
\begin{itemize}
\item[\convexity] convexity: $f_k \geq (f_{k-1}+f_{k+1})/2$ for all $k \in \{ 1,\ldots,d-2 \}$
\item[\logconvexity] logarithmic convexity: $f_k^2 \geq f_{k-1} f_{k+1}$ for all $k \in \{ 1,\ldots,d-2 \}$
\item[\unimodality] unimodality: $f_0 \leq \ldots \leq f_k \geq \ldots \geq f_{d-1}$ for some $k \in \{ 0,\ldots,d-1 \}$
\item[\barany] B\'ar\'any's property: $f_k \geq \min \{ f_0,f_{d-1} \}$ for all $k \in \{ 1,\ldots,d-2 \}$
\end{itemize}
Clearly each property implies the next one: \convexity\ $\Rightarrow$ \logconvexity\ $\Rightarrow$ \unimodality\ $\Rightarrow$ \barany.

Unimodality is known to be false in general for $d \geq 8$ and (rather trivially) true for $d \leq 4$.
For simplicial (and therefore also for simple) polytopes of arbitrary dimension a weaker version of unimodality
was proved by Bj\"orner \cite[Section 8.6]{MR96a:52011}.

Similarly, convexity is trivially true up to $d \leq 3$ and for $d=4$ follows easily from
$f_0 \geq 5$ and $f_2 \geq 2 f_3$ together with Euler's equation and duality.

\subsection*{Toric $g$-vectors} \label{sec:toric-g-vector}

To every $d$-polytope $P$ we can assign a $(\lfloor d/2 \rfloor+1)$-dimensional vector
$g(P) = (g_0^d(P), \ldots, g_{\lfloor d/2 \rfloor}^d(P))$, the \emph{toric $g$-vector} of $P$.
Its entries can be calculated via recursion \cite[Section 3.14]{MR98a:05001},
and interpreted geometrically for simplicial polytopes.
It is well known \cite{MR89f:52016} that $g_i^d(P) \geq 0$ for rational polytopes $P$
and only recently Karu \cite{MR2076929} showed that nonnegativity also holds for nonrational polytopes.

The entries of the toric $g$-vector of a polytope $P$ can be rewritten as a
linear combination of entries of the flag vector of $P$. Some special cases which
we will need are $g_0^d (P) = 1$ and $g_1^d (P) = f_0 - (d+1)$ for $d$-polytopes $P$
(note that $1 = f_\emptyset(P)$). See \cite{MR2132764} for a general description.

\subsection*{Convolutions} \label{sec:convolution}

Let $m_1$ and $m_2$ be linear forms on flag vectors of $d_1$-, resp. $d_2$-polytopes.
Then we obtain a linear form $m = m_1 * m_2$ by defining
\begin{displaymath}
  m (P) \; := \; \sum_{F \; d_1\text{-face of } P} m_1 (F) \, m_2 (P/F)
\end{displaymath}
for every $(d_1+d_2+1)$-polytope $P$.
Alternatively, the convolution can be described by defining
\begin{displaymath}
  f_S * f_T \; := \; f_{S \cup \{d_1\} \cup (T+d_1+1)}
\end{displaymath}
for $S \subseteq \{ 0,\ldots,d_1 \}$, $T \subseteq \{ 0,\ldots,d_2 \}$ (where $M+x := \{ m+x \mid m \in M \})$
and extending linearly \cite[Section 3]{MR90a:52012} \cite[Section 7]{MR2001c:52009}.
We will use this notation, occasionally writing $f_S^d$ to indicate the dimension $d$
of the polytopes the respective flag vector refers to.

\subsection*{{\cd}-index} \label{sec:cd-index}

Connected with every polytope (in fact with every Eulerian poset) is its {\cd}-index,
which is a polynomial in the non-commuting variables $\c$ and $\d$. The coefficients of the
{\cd}-index can again be viewed as linear combinations of flag vector entries \cite[Section 7]{MR2001c:52009}.
Stanley \cite{MR96b:06006} showed that the coefficients of the {\cd}-index of a polytope are nonnegative,
which again yields inequalities for the flag vector. Further useful results were obtained
by Ehrenborg \cite{MR2132764}.

From there we will adopt the following notation: write $\langle u \mid \Psi(P) \rangle$
for the coefficient of the {\cd}-monomial $u$ in the {\cd}-index of the polytope $P$.
Using linearity we can then define the number $\langle p \mid \Psi(P) \rangle$
for any {\cd}-polynomial $p$.

\bigskip

In some of the following proofs we omit the longer calculations.
For more details see the appendix.


\section{Dimension 5} \label{sec:dimension-5}

\begin{theorem} \label{thm:unimodal5}
  Unimodality \unimodality\ holds for $f$-vectors of polytopes of dimension $d \leq 5$.
\end{theorem}
\begin{proof}
  Let $P$ be a $5$-polytope and $f(P) = (f_0,f_1,f_2,f_3,f_4)$ its $f$-vector.
  Trivially, $5f_0 \leq 2f_1$ and $5f_4 \leq 2f_3$, therefore $f_0 < f_1$ and $f_3 > f_4$.
  Kalai \cite{MR90a:52012} showed that $3f_2 \geq 2f_1 + 2f_3$, hence
  \begin{displaymath}
    f_2 \; \geq \; \frac{2}{3} \, (f_1+f_3) \; > \; \frac{f_1+f_3}{2}
  \end{displaymath}
  which implies that ``there cannot be a dip'' at $f_2$.
  Therefore $f(P)$ is unimodal.
\end{proof}

\begin{theorem}
  Convexity \convexity\ fails to hold for $d \geq 5$, that is,
  the $f$-vectors of $d$-polytopes are not convex in general.
\end{theorem}
\begin{proof}
  For dimension $5$ the $f$-vector of the cyclic polytope with $n$ vertices is
  \begin{displaymath}
    f (\cyclic{5}{n}) \; = \; \left( n, \, \tfrac{n(n-1)}{2}, \, 2(n^2-6n+10), \, \tfrac{5(n-3)(n-4)}{2}, \, (n-3)(n-4) \right)
  \end{displaymath}
  (cf.~\cite[Chapter~8]{MR96a:52011}), which implies
  \begin{displaymath}
    f_1 \; = \; \frac{n^2-n}{2} \; < \; \frac{2n^2-11n+20}{2} \; = \; \frac{f_0+f_2}{2}
  \end{displaymath}
  for $n \geq 8$; see Figure \ref{fig:f-c-5-8}.
  \begin{figure}[tb]
    \centering
    \includegraphics{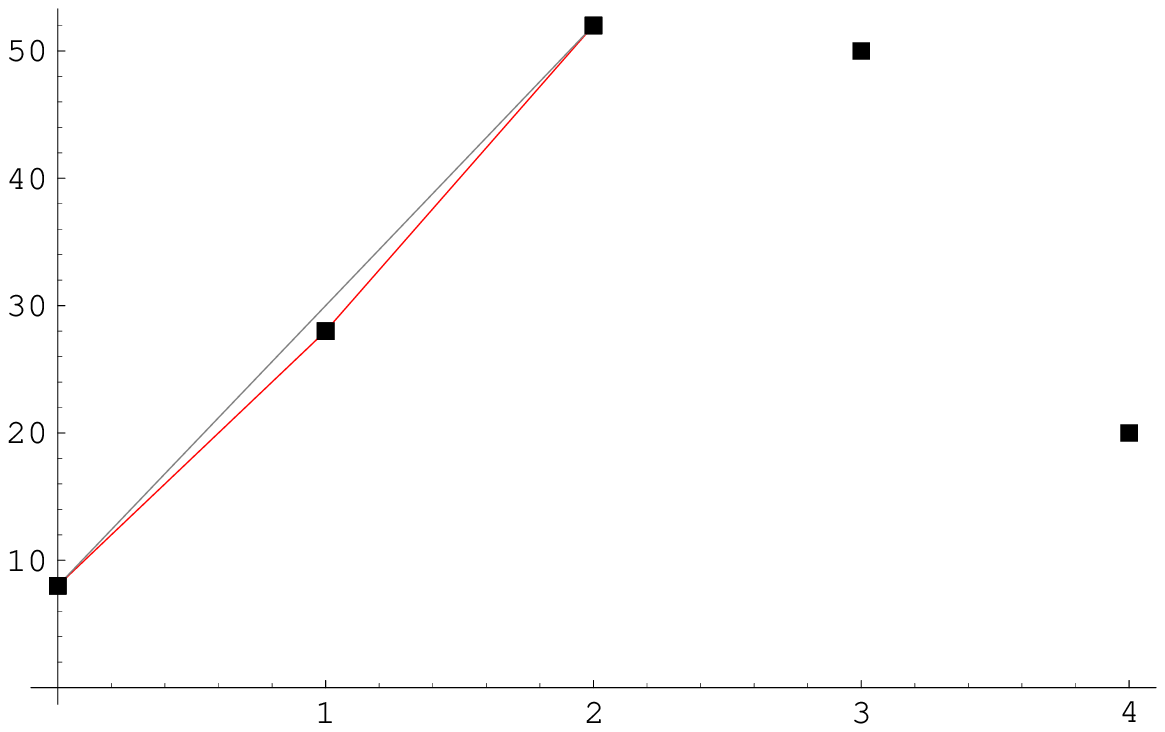}
    \caption{(non-convex) $f$-vector of $\cyclic{5}{8}$} \label{fig:f-c-5-8}
  \end{figure}
  \par
  For $d \geq 6$, cyclic $d$-polytopes are $2$-neighbourly, therefore
  $f_1 = {f_0 \choose 2}$ and $f_2 = {f_0 \choose 3}$ for $f_0 \geq d+1$. We conclude that
  \begin{displaymath}
    f_0 + f_2 - 2f_1 \; = \; \frac{1}{6} \, f_0 (f_0-2)(f_0-7) > 0
  \end{displaymath}
  for cyclic $d$-polytopes with $f_0 \geq \max\{d+2,8\}$ vertices.
  Thus for $d \geq 7$ already the simplex is a counterexample for \convexity.
\end{proof}


\section{Dimension 6} \label{sec:dimension-6}

Concerning unimodality for $f$-vectors of $6$-polytopes, we have a couple of
trivial facts, such as $f_0 < f_1$ and $f_4 > f_5$. Unimodality would therefore simply follow
from the statement $(*) \; f_1 \leq f_2$ or equivalently from $f_3 \geq f_4$ by duality.
Bj\"orner showed that the latter is true for simplicial polytopes (cf. \cite{MR96a:52011}, Theorem 8.39),
therefore in particular for cyclic polytopes, which seems to indicate that it is true in general.
However, it does not follow from the yet known inequalities -- we only have a weaker statement.

\begin{proposition} \label{prop:barany6}
  Let $f=(f_0,\ldots,f_5)$ be the $f$-vector of a $6$-polytope. Then
  \[ f_2 \; \geq \; \frac{2}{3} \, f_1 + 63 . \]
\end{proposition}
\begin{proof}
  We claim that the following inequalities hold for $f$:
  \begin{eqnarray}
    \label{in6-1} f_1 - 3f_0 & \geq & 0 \\
    \label{in6-2} f_0 - f_1 + f_2 - 21 & \geq & 0
  \end{eqnarray}
  The assertion then follows by multiplying \eqref{in6-2} by $3$ and adding \eqref{in6-1}.
  \par
  Inequality \eqref{in6-1} is trivial, simply stating that every vertex is in at least 6 edges.
  For the proof of \eqref{in6-2} we use \cite[Theorem 3.7]{MR2132764}, which implies that
  $\langle \c^2\d\c^2 - 19 \c^6 \mid \Psi(P) \rangle \geq 0$. Expressing the {\cd}-polynomial
  $\c^2\d\c^2 - 19 \c^6$ as linear combination of flag vector entries gives
  $f_0-f_1+f_2-21$ and therefore yields inequality \eqref{in6-2}.
  See the last section for detailed calculations.
\end{proof}

\begin{corollary}
  The $f$-vectors of $6$-polytopes satisfy B\'ar\'any's property \barany.
\end{corollary}
\begin{proof}
  Let $f=(f_0,\ldots,f_5)$ be the $f$-vector of a $6$-polytope.
  Clearly, $f_1 \geq 3 f_0 > f_0$, thus by Proposition~\ref{prop:barany6}
  \begin{displaymath}
    f_2 \; \geq \; \frac{2}{3} \, f_1 + 21 \; \geq \; 2 f_0 + 21 \; > \; f_0 .
  \end{displaymath}
  Dually, we have $f_3 > f_5$ and $f_4 > f_5$.
\end{proof}

As the desired inequality $(*)$ for unimodality does not follow from the known linear inequalities,
one can find vectors that satisfies all these, but not $(*)$. An example for a family of vectors is
\begin{align*}
  f^{(\ell)} = ( & f_0 , f_1 , f_2 , f_3 , f_4 ; \\
                & f_{02} , f_{03} , f_{04} , f_{13} , f_{14} , f_{24} ; \\
                & f_{024} ) \\
  \phantom{f^{(\ell)} }= ( & 22 + \ell , 111 + 3\ell , 110 + 2\ell , 35 + 4\ell , 21 + 6\ell ; \\
                & 780 + 15\ell , 1340 + 50\ell , 1080 + 51\ell , 2010 + 90\ell , 2160 + 132\ell , 1260 + 114\ell ; \\
                & 6480 + 396\ell )
\end{align*}
for $\ell \geq 0$.
The other components of these (potential) flag vectors can be calculated from the Generalized Dehn--Sommer\-ville equations.
In particular, the number of facets is $f_5 = 7+2\ell$. However it is not at all clear
that there exist polytopes having these as flag vectors.


\section{Dimension 7} \label{sec:dimension-7}

A similar statement to the one in Proposition~\ref{prop:barany6} holds for $7$-polytopes.
Nevertheless, this is not enough to prove even B\'ar\'any's property \barany, since we yet have
no condition for $f_3$.

\begin{proposition} \label{prop:barany7}
  Let $f=(f_0,\ldots,f_6)$ be the $f$-vector of a $7$-polytope. Then
  \[ f_2 \; \geq \; \frac{5}{7} \, f_1 + 36 \]
\end{proposition}
\begin{proof}
  As before, we consider two valid inequalities for $f$ which together imply the assertion:
  \begin{eqnarray}
    \label{in7-1} 2f_1 - 7f_0 & \geq & 0 \\
    \label{in7-2} f_0 - f_1 + f_2 - 36 & \geq & 0
  \end{eqnarray}
  Again, \eqref{in7-1} is trivial.
  The nonnegativity of $\langle \c^2\d\c^3-34\c^7 \mid \Psi(P) \rangle$ gives inequality \eqref{in7-2};
  see the last section.
\end{proof}

Again, one can find vectors satisfying all known linear inequalities, but violating both
$f_3 \geq f_0$ and $f_3 \geq f_6$; take, for instance, the potential flag vector
\begin{align*}
  f \; & = \; ( f_0 , f_1 , f_2 , f_3 , f_4 , f_5 ; \\
       & \qquad \quad f_{02} , f_{03} , f_{04} , f_{05} , f_{13} , f_{14} , f_{15} , f_{24} , f_{25} , f_{35} ; \\
       & \qquad \quad f_{024} , f_{025} , f_{035} , f_{135} ) \\
       & = \; ( 134 , 469 , 371 , 70 , 371 , 469 ; \\
       & \qquad \quad 2814 , 6580 , 10360 , 8484 , 9870 , 20720 , 21210 , 13790 , 20720 , 9870 ; \\
       & \qquad \quad 62160 , 84840 , 84840 , 127260 ) .
\end{align*}
From Euler's equation, we get $f_6 = 134$; nevertheless, it is again open
whether this really is the flag vector of some $7$-polytope.

As it is an open question whether logarithmic convexity holds for $f$-vectors of $7$-polytopes,
one could try to find counterexamples.
Most promising may be connected sums of cyclic polytopes,
since this construction yields counterexamples for unimodality in dimension 8
(see \cite[pp.~274f]{MR96a:52011}).

\begin{definition}
  Let $P$ and $Q$ be polytopes of the same dimension.
  If $P$ is simplicial and $Q$ simple, then a \emph{connected sum} $P \# Q$ of $P$ and $Q$ is obtained by
  cutting one vertex off $Q$ and stacking the result --- with the newly created facet --- onto $P$
  (cf.~\cite[p.~274]{MR96a:52011}).
\end{definition}

The effect of these construction on the $f$-vector of the involved polytopes
can be described as follows.

\begin{proposition} \label{prop:f-consum}
  Let $d \geq 3$ and $P$ a simplicial and $Q$ a simple $d$-polytope.
  Then the $f$-vector of $P\#Q$ is given by
  \begin{displaymath}
    f_i(P\#Q) \; = \; \left\{ \begin{array}{l@{\quad\text{ if }}l} f_i(P) + f_i(Q) & 1 \leq i \leq d-2 \\ f_i(P) + f_i(Q) - 1 & i = 0 \text{ or } i = d-1 \end{array} \right.
  \end{displaymath}
  Additionally, the $f$-vector of the connected sum $P\#\polar{P}$
  of a polytope $P$ with its dual is symmetric.
\end{proposition}
\begin{proof}
  Cutting one vertex $v$ off $Q$ decreases $f_0(Q)$ by $1$ and creates a new facet $F$,
  isomorphic to a $(d-1)$-simplex. Therefore, $f_i(Q)$ increases by ${d \choose i+1}$ if $i>0$
  and by $d-1$ if $i=0$. Afterwards all faces of both polytopes are again faces of $P\#Q$,
  except the facet $F$ in both polytopes (which completely disappears) and the new faces of $F$
  in $Q$ (which are identified with their counterparts in $P$).
  \par
  The $f$-vector of $P\#\polar{P}$ is obviously symmetric, since
  \begin{displaymath}
    f_i(P\#\polar{P}) \; = \; \left\{ \begin{array}{l@{\quad\text{ if }}l} f_i(P) + f_{d-1-i}(P) & 1 \leq i \leq d-2 \\ f_i(P) + f_{d-1-i}(P) - 1 & i = 0 \text{ or } i = d-1 \end{array} \right.
  \end{displaymath}
\end{proof}

\begin{proposition} \label{prop:logconv7}
  For all $n \geq 8$, the $f$-vector of $P_7^n := \cyclic{7}{n} \# \cyclic{7}{n}^\Delta$
  is logarithmically convex and
  \[ \frac{f_3(P_7^n)^2}{f_2(P_7^n) f_4(P_7^n)} \; \stackrel{n \rightarrow \infty}{\longrightarrow} \; 1 \]
\end{proposition}
\begin{proof}
  The proof is done by straightforward calculation; see the last section for the main steps.
\end{proof}

So in a sense, the connected sums of cyclic $7$-polytopes are as close as polytopes can get to
logarithmic non-convexity.


\section{Summary} \label{sec:summary}

The results can be summarized as in Table~\ref{tab:summary}.
\begin{table}[bt]
  \centering
  \begin{tabular}{|c|ccccc|} \hline
    Dimension     & $\leq 4$ &       $5$ &       $6$ &       $7$ &  $\geq 8$ \\ \hline
    \convexity    & \holds   & \holdsnot & \holdsnot & \holdsnot & \holdsnot \\
    \logconvexity & \holds   & ?         & ?         & ?         & \holdsnot \\
    \unimodality  & \holds   & \holds    & ?         & ?         & \holdsnot \\
    \barany       & \holds   & \holds    & \holds    & ?         & ?         \\ \hline
  \end{tabular}
  \caption{Summary of known properties for polytopes --- a \holds, resp.~\holdsnot\/ indicates
    that the given property holds, resp.~does not hold for all polytopes of the given dimension.}
  \label{tab:summary}
\end{table}
In the light of Proposition~\ref{prop:logconv7}, the following conjecture seems natural.

\begin{conjecture}
  \logconvexity\ holds for $d$-polytopes of dimension $d \leq 7$.
\end{conjecture}


\begin{appendix}

\section*{Detailed calculations} \label{sec:deta-calc}

\subsection*{Proof of Theorem \ref{thm:unimodal5} \\ {\normalsize \rm (cf.~\cite[Section 7]{MR90a:52012})}} \label{sec:proof-theor-refthm}

Let $P$ be a $5$-polytope.
\begin{eqnarray*}
  (g_0^1 * g_1^2 * g_0^0)(P) & = & f_\emptyset^1 * (f_0^2-3) * f_\emptyset^0 \; = \; (f_{12}^4 - 3f_1^4) * f_\emptyset^0 \; = \; f_{124} - 3f_{14} \\
    & = & f_{123} - 3 \, (2f_1 - f_{12} + f_{13}) \; = \; -6f_1 + 3f_{02} - f_{13} \\
  (g_0^0 * g_1^2 * g_0^1)(P) & = & f_\emptyset^0 * (f_0^2-3) * f_\emptyset^1 \; = \; (f_{01}^3 - 3f_0^3) * f_\emptyset^1 \; = \; f_{013} - 3f_{03} \\
    & = & 2f_{13} - 3f_{03} \\
  (g_1^2 * g_1^2)(P) & = & (f_0^2 - 3) * (f_0^2 - 3) \; = \; f_{023} - 3f_{02} - 3f_{23} + 9f_2 \\
    & = & f_{013} - 3f_{02} - 3 \, (2f_3 - f_{03} + f_{13}) + 9f_2 \\
    & = & -f_{13} - 3f_{02} - 6f_3 + 3f_{03} + 9f_2
\end{eqnarray*}
by the rules of convolution and the Generalized Dehn--Sommerville equations \cite{MR86f:52010b}.
Hence we have
\begin{displaymath}
  \begin{array}{rrcrcrcrcrcl}
    -6 f_1 &       &   &       & + & 3 f_{02} &   &          & - & f_{13}   & \geq & 0 \\
           &       &   &       &   &          & - & 3 f_{03} & + & 2 f_{13} & \geq & 0 \\
           & 9 f_2 & - & 6 f_3 & - & 3 f_{02} & + & 3 f_{03} & - & f_{13}   & \geq & 0
  \end{array}
\end{displaymath}
Adding up all three inequalities yields $ -6 f_1 + 9 f_2 - 6 f_3 \geq 0 $,
that is the assertion $3f_2 \geq 2f_1 + 2f_3$. \qed

\subsection*{Proof of Inequality~\eqref{in6-2}} \label{sec:proof:prop:barany6}

We express the {\cd}-word $\c^2\d\c^2$ in terms of the flag vector of the $6$-polytope $P$
by applying \cite[Proposition~7.1]{MR2001c:52009}:
\begin{displaymath}
  \langle \c^2\d\c^2 \mid \Psi(P) \rangle \; = \; \sum_{i=0}^2 (-1)^{4-i} k_i \; = \; k_0 - k_1 + k_2 .
\end{displaymath}
For the sparse flag $k$-vector we have
\begin{displaymath}
  k_i \; = \; \sum_{T \subseteq \{i\}} (-2)^{1-|T|} f_T \; = \; -2 f_\emptyset + f_i \; = \; f_i - 2
\end{displaymath}
and therefore $\langle \c^2\d\c^2 \mid \Psi(P) \rangle = f_0-f_1+f_2-2$. The trivial {\cd}-word
$\c^6$ translates into $f_\emptyset = 1$, hence
\begin{displaymath}
  \langle \c^2\d\c^2 - 19\c^6 \mid \Psi(P) \rangle \; = \; f_0-f_1+f_2-21 .
\end{displaymath} \qed

\subsection*{Proof of Inequality~\eqref{in7-2}} \label{sec:proof:prop:barany7}

We calculate for the $7$-polytope $P$ exactly as above (the additional $\c$ at the end
has no influence whatsoever on the calculation):
\begin{displaymath}
  \langle \c^2\d\c^3 \mid \Psi(P) \rangle \; = \; f_0-f_1+f_2-2 .
\end{displaymath}
Together with $\c^7$, which again represents $f_\emptyset$, we get
\begin{displaymath}
  \langle \c^2\d\c^3 - 34\c^7 \mid \Psi(P) \rangle \; = \; f_0-f_1+f_2-36 .
\end{displaymath} \qed

\subsection*{Proof of Proposition~\ref{prop:logconv7}}

The $f$-vector of the cyclic $7$-polytope on $n$ vertices is given by
\begin{align*}
  f(\cyclic{7}{n}) \; = \; \big( n & , \, \tfrac{n(n-1)}{2} , \, \tfrac{n(n-1)(n-2)}{6} , \, \tfrac{5(n-4)(n^2-8n+21)}{6} , \\
  & \tfrac{(n-4)(3n^2-31n+84)}{2} , \, \tfrac{7(n-4)(n-5)(n-6)}{6} , \, \tfrac{(n-4)(n-5)(n-6)}{3} \big)
\end{align*}
(cf.~\cite[Chapter~8]{MR96a:52011}).
From this we obtain $f(P_7^n)=(f_0(n),\ldots,f_6(n))$ by Proposition~\ref{prop:f-consum}:
\begin{align*}
  f_0(n) \; = \; \tfrac{(n-3)(n^2-12n+41)}{3} \quad , & \quad f_1(n) \; = \; \tfrac{7n^3-102n^2+515n-840}{6} , \\
  f_2(n) \; = \; \tfrac{5n^3-66n^2+313n-504}{3} \quad , & \quad f_3(n) \; = \; \tfrac{5(n-4)(n^2-8n+21)}{3}
\end{align*}
By symmetry of $f(P_7^n)$, these entries suffice to verify logarithmic convexity.
We get
\begin{align*}
  \frac{f_1(n)^2}{f_0(n)f_2(n)} \; & = \; \frac{(7n^3-102n^2+515n-840)^2}{4(n-3)(n^2-12n+41)(5n^3-66n^2+313n-504)} \; > \; 1 \, , \\
  \frac{f_2(n)^2}{f_1(n)f_3(n)} \; & = \; \frac{2(5n^3-66n^2+313n-504)^2}{5(n-4)(n^2-8n+21)(7n^3-102n^2+515n-840)} \; > \; 1 \, , \\
  \frac{f_3(n)^2}{f_2(n)f_4(n)} \; & = \; \frac{25(n-4)^2(n^2-8n+21)^2}{(5n^3-66n^2+313n-504)^2} \; > \; 1
\end{align*}
for $n \geq 8$. Since the leading coefficients of the polynomials in the numerator and the denominator
of the last fraction are equal,
\begin{displaymath}
  \frac{f_3(P_7^n)^2}{f_2(P_7^n) f_4(P_7^n)} \; \stackrel{n \rightarrow \infty}{\longrightarrow} \; 1 .
\end{displaymath} \qed

\end{appendix}

\vfill \bibliographystyle{siam}

\end{document}